\documentclass[preprint,12pt]{elsarticle}
\journal{J.}
\usepackage{dsfont}
\usepackage{amsfonts}
\usepackage{amsmath}
\usepackage{epsfig}
\usepackage{amssymb}
\begin{document}
\bibliographystyle{elsarticle-harv}
\begin{frontmatter}
\title{Stochastic minimum-energy control}
\author{Bujar Gashi\footnote{Department of Mathematical Sciences, The University of Liverpool, Liverpool, L69 7ZL, UK; Email: Bujar.Gashi@liverpool.ac.uk}}

\begin{abstract}
We give the solution to the minimum-energy control problem for linear stochastic systems. The problem is as follows: given an exactly controllable system, find the control process with the minimum expected energy that transfers the system from a given initial state to a desired final state. The solution is found in terms of a certain  forward-backward stochastic differential equation of Hamiltonian type.
\end{abstract}


\end{frontmatter}

\numberwithin{equation}{section}
\newtheorem{proof}{Proof}
\newtheorem{definition}{Definition}
\newtheorem{theorem}{Theorem}
\newtheorem{lemma}{Lemma}
\newtheorem{proposition}{Proposition}
\newtheorem{remark}{Remark}
\newtheorem{corollary}{Corollary}
\newtheorem{assumption}{Asssumption}
\section{Introduction} The notion of {\it controllability} was introduced by Kalman~\cite{kalman1} and it characterizes the ability of controls to transfer a system from a given initial state to a desired final state. When the system is completely controllable, there are many controls that can achieve such a transfer of the system state.  This naturally leads to the problem of choosing the ``best'' control that performs this task. Another contribution of Kalman~\cite{kalman1} was the solution of this problem for linear deterministic systems and using the quadratic cost as an optimality criterion. These kinds of problems are know as the {\it minimum-energy control} problems (see also~\cite{kalman2},~\cite{kalman3},~\cite{athans},~\cite{lewis}).\\

In this paper we formulate and solve the minimum-energy control problem for linear {\it stochastic} systems driven by a Brownian motion. The notion of controllability that we adapt is that of {\it exact controllability}, as introduced by Peng~\cite{peng} (see also~\cite{liu1},~\cite{liu2},~\cite{GP},~\cite{wang}). This notion of controllability is a faithful extension of Kalman's notion of complete controllability to stochastic systems. The difference between these two definitions is that in the case of exact controllability the terminal state can be a random variable rather than a fixed number. This makes the stochastic minimum-energy control problem considerably harder than in deterministic setting.\\

The precise formulation of the stochastic minimum-energy control problem is given in the next section. This is followed by the proof of solvability for a Hamiltonian system and its relation with exact controllability. Section \ref{minsec} contains the solution to the stochastic minimum energy control problem. As an extension of this result, we give the solution to the stochastic linear-quadratic (LQ) regulator problem with a fixed final state in the final section.

\section{Problem formulation}
Let $(\Omega, \mathcal{F}, (\mathcal{F}_t,t\geq0),\mathds{P})$ be a given complete filtered probability space on which the scalar standard Brownian motion $(W(t), t\geq 0)$ is defined. We assume that $\mathcal{F}_t$ is the augmentation of $\sigma\{W(s):0\leq s\leq t\}$ by all the $\mathds{P}$-null sets of $\mathcal{F}$. If $\xi:\Omega\rightarrow\mathds{R}^n$ is an $\mathcal{F}_T$-measurable random variable such that $\mathds{E}[|\xi|^2]<\infty$, we write $\xi\in L^2(\Omega, \mathcal{F}_T,\mathds{P}; \mathds{R}^n)$. If $f:[0,T]\times\Omega\rightarrow\mathds{R}^n$ is an $\{\mathcal{F}_t\}_{t\geq 0}$ adapted process and if $\mathds{E}\int_0^T|f(t)|^2dt<\infty$, we write $f(\cdot)\in L_{\mathcal{F}}^2(0,T;\mathds{R}^n)$; if $f(\cdot)$ has a.s. continuous sample paths and $\mathds{E}\sup_{t\in[0,T]}|f(t)|^2<\infty$, we write $f(\cdot)\in L_{\mathcal{F}}^2(\Omega;C(0,T;\mathds{R}^n))$; if $f(\cdot)$ is uniformly bounded (i.e. esssup$_{t\in[0,T]}|f(t)|$ $<\infty$), we write $f(\cdot)\in L^\infty(0,T;\mathds{R}^n)$.\\

Consider the linear stochastic control system:
\begin{eqnarray}
\left\{
\begin{array}{l}
\displaystyle dx(t)=[A(t)x(t)+B(t)u(t)]dt+[C(t)x(t)+D(t)u(t)]dW(t)\\
\\
\displaystyle x(0)=x_0\in\mathds{R}^n,\quad\mbox{is given}.
\end{array}\right.\label{main system}
\end{eqnarray}
We assume that $A(\cdot), C(\cdot)\in L^\infty(0,T;\mathds{R}^{n\times n})$, and $B(\cdot), D(\cdot)\in L^\infty (0,T;\mathds{R}^{n\times m})$.  If the control process $u(\cdot)$ belongs to $L_{\mathcal{F}}^2(0,T;\mathds{R}^m)$, then \eqref{main system} has a unique strong solution $x(\cdot)\in L_{\mathcal{F}}^2(\Omega;C(0,T;\mathds{R}^n))$ (see, e.g. Theorem 1.6.14 of~\cite{yong}).\\

For a given $\xi\in L^2(\Omega,\mathcal{F}_T,P;\mathds{R}^n)$, we are interested in the following subset of control processes:
\begin{eqnarray}
\mathcal{U}_{\xi}\equiv\left\{u(\cdot)\in L^2_\mathcal{F}(0,T;\mathds{R}^m):x(T)=\xi\quad a.s.\right\}.\nonumber
\end{eqnarray}

{\bf Minimum-energy control problem.} {\it Let $R(\cdot)\in L^\infty(0,T;\mathds{R}^{m\times m})$ be a given symmetric matrix such that $R(t)>0$, $a.e.$ $t\in[0,T]$. For any given $x_0\in\mathds{R}^n$ and $\xi\in L^2(\Omega,\mathcal{F}_T,P;\mathds{R}^n)$ find the control process $u(\cdot)\in\mathcal{U}_{\xi}$ that minimizes the cost functional}
\begin{eqnarray}
J(u(\cdot))=\mathds{E}\int_0^Tu'(t)R(t)u(t)dt.\label{main cost}
\end{eqnarray}

This is clearly the stochastic version of the Kalman's minimum energy control problem. A related problem was considered by Klamka~\cite{klamka}. However, in~\cite{klamka} only the linear stochastic systems  with {\it additive} noise are considered, whereas (\ref{main system}) has a {\it multiplicative} noise. Our approach to solving the stochastic minimum-energy control problem is different from that of~\cite{klamka}, where an operator-theoretic method was used, whereas here we base our approach on a forward-backward stochastic differential equation of a Hamiltonian type.\\

In order to ensure that the set $\mathcal{U}_{\xi}$ is not empty, we make some assumptions on the controllability of (\ref{main system}). Out of the many possible notions of controllability for stochastic systems, we employ the notion of {\it exact controllability} as introduced by Peng~\cite{peng}.

\begin{definition} System (\ref{main system}) is called exactly controllable at time $T>0$ if for any $x_0\in\mathds{R}^n$ and $\xi\in L^2(\Omega,\mathcal{F}_T,\mathds{P};\mathds{R}^n)$, there exists at least one control $u(\cdot)\in L_{\mathcal{F}}^2(0,T;\mathds{R}^m)$, such that the corresponding trajectory $x(\cdot)$ satisfies the initial condition $x(0)=x_0$ and the terminal condition $x(T)=\xi$, $a.s.$.
\end{definition}

We solve the minimum-energy control problem under the following two assumptions.\\
\\
{\bf (A1)} The system \eqref{main system} is exactly controllable at time $T>0$.\\
\\
{\bf (A2)} There exists an {\it invertible} matrix $M(\cdot)\in L^\infty(0,T;\mathds{R}^{m\times m})$ such that $D(t)M(t)=[I,0]$.\\
\\
Assumption A1 ensures that the set $\mathcal{U}_{\xi}$ is not empty. Assumption A2 implies that $m\geq n$, i.\ e. the number of control inputs to the system is at least as large as the number of the states of the system. This may appear as a strong assumption when compared with the minimum-energy control problem of deterministic systems. However, at least when the matrix $D(\cdot)$ has continuous coefficients, this assumption is implied by assumption A1. Indeed, by Proposition 2.1. of~\cite{peng}, a {\it necessary} condition for exact controllability at time $T$ of the system (\ref{main system}) is that $rank$ $D(t)=n$, $\forall t\in[0,T]$. Then from the Dole\v{z}al's theorem~\cite{dolezal}, it follows that there exists the matrix $M(\cdot)$ in assumption A2.\\

We now reformulate the minimum-energy control problem in a more convenient form. Let the processes $z(\cdot)\in L_{\mathcal{F}}^2(0,T;\mathds{R}^n)$ and $v(\cdot)\in L_{\mathcal{F}}^2(0,T;\mathds{R}^{m-n})$ be such that
\begin{eqnarray}
u(t)=M(t)\left[
\begin{array}{l}
z(t)\\
\\
v(t)
\end{array}\right].\label{uvz}
\end{eqnarray}
Let the matrices $G(\cdot)\in L^\infty(0,T;\mathds{R}^{n\times n})$, $F(\cdot)\in L^\infty(0,T;\mathds{R}^{n\times (m-n)})$,  $H_1(\cdot)\in L^\infty(0,T;\mathds{R}^{n\times n})$, $H_2(\cdot)\in L^\infty(0,T;\mathds{R}^{n\times (m-n)})$, $H_3(\cdot)\in L^\infty(0,T;\mathds{R}^{(m-n)^2})$, be such that
\begin{eqnarray}
B(t)M(t)=\left[
\begin{array}{ll}
G(t)&F(t)
\end{array}\right],\quad M'(t)R(t)M(t)=\left[
\begin{array}{ll}
H_1(t)&H_2(t)\\
\\
H_2'(t)&H_3(t)
\end{array}\right].\label{decomp}
\end{eqnarray}
Due to the symmetric nature of the matrix $R(\cdot)$, the matrices $H_1(\cdot)$ and $H_3(\cdot)$ are also symmetric. Moreover, due to the positive definiteness of $R(\cdot)$ and the Schur's lemma, it holds that
\begin{eqnarray}
&&H_3(t)>0,\quad a.e.\quad t\in[0,T],\nonumber\\
\nonumber\\
&&H_1(t)-H_2'(t)H_3^{-1}(t)H_2(t)>0,\quad a.e.\quad t\in[0,T].\nonumber
\end{eqnarray}
Equation (\ref{main system}) and the cost functional (\ref{main cost}) can now be written as
\begin{eqnarray}
\left\{
\begin{array}{l}
\displaystyle dx(t)=[A(t)x(t)+F(t)v(t)+G(t)z(t)]dt+[C(t)x(t)+z(t)]dW(t),\\
\\
\displaystyle x(0)=x_0\in\mathds{R}^n,\quad\mbox{is given},
\end{array}\right.\label{main2}
\end{eqnarray}
\begin{eqnarray}
J(v(\cdot),z(\cdot))\!\!=\!\!\mathds{E}\!\!\int_0^T\!\![z'(t)H_1(t)z(t)\!+\!2v'(t)H_2'(t)z(t)\!+\!v'(t)H_3(t)v(t)]dt.\label{c}
\end{eqnarray}
To each element of the set $\mathcal{U}_\xi$ it corresponds a pair of processes $(v(\cdot),z(\cdot))$ from the set
\begin{eqnarray}
\mathcal{A}_\xi\equiv\left\{v(\cdot)\in L_{\mathcal{F}}^2(0,T;\mathds{R}^{m-n}), z(\cdot)\in L_{\mathcal{F}}^2(0,T;\mathds{R}^n): x(T)=\xi\quad a.s. \right\}.\nonumber
\end{eqnarray}
In this reformulation, the minimum-energy control problem is:
\begin{eqnarray}
\left\{
\begin{array}{l}
\displaystyle \min_{(v(\cdot),z(\cdot))\in\mathcal{A}_\xi}J(v(\cdot),z(\cdot)),\\
\\
\displaystyle s.t.\quad\mbox{(\ref{main2})}.
\end{array}\right.\label{ME}
\end{eqnarray}
Before we proceed to its solution, let us state a useful necessary and sufficient condition for the exact controllability of (\ref{main2}). It is a slight modification of the result in~\cite{liu1}, and we thus omit the proof.
\begin{proposition}Let $E(\cdot)\in L^\infty(0,T;\mathds{R}^{m\times m})$ be any symmetric matrix such that $E(t)>0$, $a.e.$ $t\in[0,T]$. Also let $\Phi(\cdot)$ be the unique solution to the equation
\begin{eqnarray}
\left\{
\begin{array}{l}
\displaystyle d\Phi(t)=-\Phi(t)[A(t)-G(t)C(t)]dt-\Phi(t)G(t)dW(t),\\
\\
\displaystyle \Phi(0)=I.
\end{array}\right.\nonumber
\end{eqnarray}
The system (\ref{main2}) is exactly controllable at time $T$ if and only if
\begin{eqnarray}
rank\quad\left[\mathds{E}\int_0^T\Phi(t)F(t)E(t)F'(t)\Phi'(t)dt\right]=n.\label{ra}
\end{eqnarray}
\end{proposition}
\section{Hamiltonian system of equations}\label{hamsec}
The following forward-backward stochastic differential equation  of Hamiltonian type appears naturally in the next section:
\begin{eqnarray}
\left\{
\begin{array}{l}
\displaystyle dX(t)=\{A(t)X(t)+F(t)H_3^{-1}(t)F'(t)Y(t)+[G(t)-F(t)H_3^{-1}(t)H_2'(t)]Z(t)\}dt\\
\\
\displaystyle \quad\quad\quad+[C(t)X(t)+Z(t)]dW(t),\\
\\
\displaystyle dY(t)=-[A'(t)-C'(t)G'(t)+C'(t)H_2(t)H_3^{-1}(t)F'(t)]Y(t)dt\\
\\
\displaystyle \quad\quad\quad -C'(t)[H_1(t)-H_2(t)H_3^{-1}(t)H_2'(t)]Z(t)dt\\
\\
\displaystyle\quad\quad\quad +\{[-G'(t)+H_2(t)H_3^{-1}(t)F'(t)]Y(t)+[H_1(t)-H_2(t)H_3^{-1}(t)H_2'(t)]Z(t)\}dW(t),\\
\\
\displaystyle X(0)=x_0,\quad X(T)=\xi,\quad Y(0)=K,
\end{array}\right.\label{hamiltonian}
\end{eqnarray}
where we can choose the vector $K\in\mathds{R}^n$. To simplify the notation, we introduce the matrices:
\begin{eqnarray}
\bar{A}(t)&\equiv& A(t)-G(t)C(t)+F(t)H_3^{-1}(t)H_2'(t)C(t),\nonumber\\
\nonumber\\
\bar{B}(t)&\equiv& G(t)-F(t)H_3^{-1}(t)H_2'(t),\nonumber\\
\nonumber\\
\bar{H}(t)&\equiv& H_1(t)-H_2(t)H_3^{-1}(t)H_2'(t).\nonumber
\end{eqnarray}
By defining $\bar{X}(t)\equiv-X(t)$ and $\bar{Z}(t)\equiv-[C(t)X(t)+Z(t)]$, we can rewrite (\ref{hamiltonian}) as
\begin{eqnarray}
\left\{
\begin{array}{l}
\displaystyle d\bar{X}(t)=[\bar{A}(t)\bar{X}(t)-F(t)H_3^{-1}(t)F'(t)Y(t)+\bar{B}(t)\bar{Z}(t)]dt+\bar{Z}(t)dW(t),\\
\\
\displaystyle dY(t)=[-\bar{A}'(t)Y(t)-C'(t)\bar{H}(t)C(t)\bar{X}(t)+C'(t)\bar{H}(t)\bar{Z}(t)]dt\\
\\
\displaystyle\quad\quad\quad +[-\bar{B}'(t)Y(t)+\bar{H}(t)C(t)\bar{X}(t)-\bar{H}(t)\bar{Z}(t)]dW(t),\\
\\
\displaystyle \bar{X}(0)=-x_0,\quad \bar{X}(T)=-\xi,\quad Y(0)=K.
\end{array}\right.\label{hamiltonian2}
\end{eqnarray}
This forward-backward stochastic differential equation is similar to the Hamiltonian system of stochastic LQ control problem~\cite{yong}. Two main differences are that here the initial value of of $\bar{X}(\cdot)$ is fixed and the vector $K$ can be chosen. We thus seek the solution {\it quadruple} $(\bar{X}(t),Y(t),\bar{Z}(t),K)$, rather than the solution triple $(\bar{X}(t),Y(t),\bar{Z}(t))$, as is usually the case with Hamiltonian systems.
\begin{theorem} There exists a unique solution quadruple $(\bar{X}(\cdot),Y(\cdot),\bar{Z}(\cdot),K)\in L_\mathcal{F}^2(0,T;\mathds{R}^n)\times L_\mathcal{F}^2(0,T;\mathds{R}^n)\times L_\mathcal{F}^2(0,T;\mathds{R}^n)\times\mathds{R}^n$ to \eqref{hamiltonian2} for any $x_0\in\mathds{R}^n$ and any $\xi\in L^2_\mathcal{F}(\Omega,\mathds{P},\mathcal{F};\mathds{R}^n)$, if and only if the system \eqref{main2} is exactly controllable at time $T$. In this case, $X(\cdot)$ and $Z(\cdot)$ in terms of $Y(\cdot)$, and the explicit formula for $K$, are given by
\begin{eqnarray}
\bar{X}(t)&=&\bar{P}(t)Y(t)+p(t),\label{X}\\
\nonumber\\
\bar{Z}(t)&=&[I+\bar{P}(t)\bar{H}(t)]^{-1}[\bar{P}(t)\bar{H}(t)C(t)\bar{P}(t)-\bar{P}(t)\bar{B}'(t)]Y(t)\nonumber\\
\nonumber\\
&+&[I+\bar{P}(t)\bar{H}(t)]^{-1}[\bar{P}(t)\bar{H}(t)C(t)p(t)+q(t)],\label{Z}\\
\nonumber\\
K&=&\bar{P}^{-1}(0)\left\{x_0-\mathds{E}[\mathcal{P}(T)\xi]\right\}.\label{K}
\end{eqnarray}
Here $\bar{P}(\cdot)$ is the unique solution of the Riccati equation
\begin{eqnarray}
\left\{
\begin{array}{l}
\displaystyle \dot{\bar{P}}(t)-A(t)\bar{P}(t)-\bar{P}(t)A'(t)+F(t)H_3^{-1}(t)F'(t)+\bar{B}(t)\bar{H}^{-1}(t)\bar{B}'(t)\\
\\
\displaystyle -[\bar{P}(t)C'(t)-\bar{B}\bar{H}^{-1}][\bar{H}^{-1}(t)+\bar{P}(t)]^{-1}[\bar{P}(t)C'(t)-\bar{B}\bar{H}^{-1}]'=0,\\
\\
\displaystyle \bar{P}(T)=0,
\end{array}\right.\label{RiccatiH}
\end{eqnarray}
whereas $(p(\cdot),q(\cdot))$ are the unique solution pair of the following linear backward stochastic differential equation
\begin{eqnarray}
\left\{
\begin{array}{l}
\displaystyle dp(t)=[\mathcal{B}_1(t)p(t)+\mathcal{B}_2(t)q(t)]dt+q(t)dW(t),\\
\\
\displaystyle p(T)=\xi,\\
\\
\displaystyle \mathcal{B}_1(t)\equiv\bar{A}(t)+\bar{P}(t)C'(t)\bar{H}(t)C(t)\\
\\
\displaystyle\quad\quad\quad+[\bar{B}(t)-\bar{P}(t)C'(t)\bar{H}(t)][I+\bar{P}(t)\bar{H}(t)]^{-1}\bar{P}(t)\bar{H}(t)C(t),\\
\\
\displaystyle \mathcal{B}_2(t)\equiv[\bar{B}(t)-\bar{P}(t)C'(t)\bar{H}(t)][I+\bar{P}(t)\bar{H}(t)]^{-1}.
\end{array}\right.\label{bsde}
\end{eqnarray}
Finally, the process $\mathcal{P}(\cdot)$ is the unique solution to the stochastic differential equation
\begin{eqnarray}
\left\{
\begin{array}{l}
\displaystyle d\mathcal{P}(t)=-\mathcal{P}(t)\mathcal{B}_1(t)dt-\mathcal{P}(t)\mathcal{B}_2(t)dW(t),\\
\\
\displaystyle \mathcal{P}(0)=I.
\end{array}\right.\label{mathcalP}
\end{eqnarray}
\end{theorem}
 {\bf Proof.} ({\it Positivity of $\bar{P}(0)$}) The Riccati differential equation \eqref{RiccatiH} has a unique solution $\bar{P}(t)\geq0$, $\forall t\in[0,T]$ (see, e.g., Theorem 3.1 of~\cite{XYZ}). We show that $\bar{P}(0)>0$ if and only if \eqref{main2} is exactly controllable at time $T$. Thus consider the stochastic control system
\begin{eqnarray}
\left\{
\begin{array}{l}
\displaystyle d\widetilde{x}(t)=[-A'(t)\widetilde{x}(t)+C'(t)\widetilde{u}(t)]dt-\widetilde{u}(t)dW(t),\\
\\
\displaystyle \widetilde{x}(0)=\widetilde{x}_0\neq0,
\end{array}\right.\label{tx}
\end{eqnarray}
and the associated cost functional
\begin{eqnarray}
\widetilde{J}(\widetilde{u}(\cdot))&=&\mathds{E}\int_0^T\widetilde{x}'(t)[F(t)H_3^{-1}(t)F'(t)+\bar{B}(t)\bar{H}^{-1}(t)\bar{B}'(t)]\widetilde{x}(t)dt\nonumber\\
\nonumber\\
&+&\mathds{E}\int_0^T[-2\widetilde{x}'(t)\bar{B}(t)\bar{H}^{-1}(t)\widetilde{u}(t)+\widetilde{u}'(t)\bar{H}^{-1}(t)\widetilde{u}(t)]dt.\label{tc}
\end{eqnarray}
From Theorem 2.2 of~\cite{XYZ}, it follows that $\widetilde{x}'_0\bar{P}(0)\widetilde{x}_0=\min_{\tilde{u}(\cdot)} \widetilde{J}(\widetilde{u}(\cdot))$. Introducing the new control $\tilde{v}(t)\equiv\tilde{u}(t)-\bar{B}'\tilde{x}(t)$ transforms \eqref{tx} and \eqref{tc} into
\begin{eqnarray}
\left\{
\begin{array}{l}
\displaystyle d\widetilde{x}(t)=[-A'(t)+C'(t)\bar{B}'(t)]\widetilde{x}(t)dt+C'(t)\widetilde{v}(t)dt-[\bar{B}'\tilde{x}(t)+\widetilde{v}(t)]dW(t),\\
\\
\displaystyle \widetilde{x}(0)=\widetilde{x}_0\neq0,
\end{array}\right.\label{tx2}
\end{eqnarray}
\begin{eqnarray}
\widetilde{J}(\widetilde{v}(\cdot))=\mathds{E}\int_0^T[\widetilde{x}'(t)F(t)H_3^{-1}(t)F'(t)\widetilde{x}(t)+\widetilde{v}'(t)\bar{H}^{-1}(t)v(t)]dt.\label{tc2}
\end{eqnarray}

To prove the {\it sufficiency}, let us assume the opposite, i.e. that $\widetilde{x}'_0\bar{P}(0)\widetilde{x}_0=0$ and the system \eqref{main2} is exactly controllable at time $T$. From \eqref{tc2}, and the fact that $\bar{H}^{-1}(t)>0$, $a.e.$ $t\in[0,T]$, we conclude that in order to minimize (\ref{tc2}) it is necessary to have $\widetilde{v}(t)=0$, $a.e.$ $t\in[0,T]$ $a.s.$. For such a $\widetilde{v}(\cdot)$, the solution to \eqref{tx2} becomes $\widetilde{x}(t)=\widetilde{\Phi}(t)x_0$, with $\widetilde{\Phi}(t)$ being the solution to the stochastic differential equation
\begin{eqnarray}
\left\{
\begin{array}{l}
\displaystyle d\widetilde{\Phi}(t)=[-A'(t)+C'(t)\bar{B}'(t)]\widetilde{\Phi}(t)dt-\bar{B}'\widetilde{\Phi}(t)dW(t),\\
\\
\displaystyle \widetilde{\Phi}(0)=I,
\end{array}\right.\label{Pt}
\end{eqnarray}
The cost \eqref{tc2} now becomes
\begin{eqnarray}
0=\widetilde{x}_0'\left[\mathds{E}\int_0^T\widetilde{\Phi}'(t)F(t)H_3^{-1}(t)F'(t)\widetilde{\Phi}(t)dt\right]\widetilde{x}_0,\label{prt}
\end{eqnarray}
which means condition (\ref{ra}) does not hold. Since (\ref{ra}) is necessary for the exact controllability of \eqref{main2}, we have a contradiction.\\

To prove the {\it necessity}, let us assume that $\bar{P}(0)>0$ and the system \eqref{main2} is not exactly controllable at time $T$. Taking $\widetilde{v}(t)=0$ $a.e.$ $t\in[0,T]$ $a.s.$, the cost \eqref{tc2} becomes the right-hand side of \eqref{prt}, and thus there exists $\widetilde{x}_0\neq 0$ such that equation \eqref{prt} holds. This means that $\min_{\tilde{v}(\cdot)} \widetilde{J}(\widetilde{v}(\cdot))=0$, which contradicts the fact that is should be $\widetilde{x}'_0\bar{P}(0)\widetilde{x}_0>0$.\\

({\it Solvability of (\ref{hamiltonian2})}) Here we follow Yong~\cite{yong4},~\cite{yong3}, in seeking the relation
\begin{eqnarray}
\bar{X}(t)&=&\bar{P}(t)Y(t)+p(t),\label{spt}
\end{eqnarray}
which implies that the differential of $\bar{X}$ should be
\begin{eqnarray}
&&d\bar{X}(t)=\dot{\bar{P}}(t)Y(t)+\bar{P}(t)dY(t)+dp(t)\nonumber\\
\nonumber\\
&&=\dot{\bar{P}}(t)Y(t)dt+\bar{P}(t)[-\bar{A}'(t)Y(t)-C'(t)\bar{H}(t)C(t)\bar{X}(t)+C'(t)\bar{H}(t)\bar{Z}(t)]dt\nonumber\\
\nonumber\\
&&+[\mathcal{B}_1(t)p(t)+\mathcal{B}_2(t)q(t)]dt\nonumber\\
\nonumber\\
&&+\bar{P}(t)[-\bar{B}'(t)Y(t)+\bar{H}(t)C(t)\bar{X}(t)-\bar{H}(t)\bar{Z}(t)]dW(t)+q(t)dW(t)\nonumber
\end{eqnarray}
By comparing this differential with $d\bar{X}(t)$ in (\ref{hamiltonian2}), we conclude that for $a.e.$ $t\in[0,T]$ $a.s.$, we must have
\begin{eqnarray}
\left\{
\begin{array}{l}
\displaystyle \bar{A}(t)\bar{X}(t)-F(t)H_3^{-1}(t)F'(t)Y(t)+\bar{B}(t)\bar{Z}(t)=\mathcal{B}_1(t)p(t)+\mathcal{B}_2(t)q(t)\\
\\
\displaystyle +\dot{\bar{P}}(t)Y(t)+\bar{P}(t)[-\bar{A}'(t)Y(t)-C'(t)\bar{H}(t)C(t)\bar{X}(t)+C'(t)\bar{H}(t)\bar{Z}(t)],
\end{array}\right.\label{Xeq}
\end{eqnarray}
\begin{eqnarray}
\bar{Z}(t)=\bar{P}(t)[-\bar{B}'(t)Y(t)+\bar{H}(t)C(t)\bar{X}(t)-\bar{H}(t)\bar{Z}(t)]+q(t).\label{Zeq}
\end{eqnarray}
Since the matrix $[I+\bar{P}(t)\bar{H}(t)]=[\bar{H}^{-1}(t)+\bar{P}(t)]\bar{H}(t)$ is invertible, from \eqref{Zeq} we obtain \eqref{Z}. Substituting such a $\bar{Z}(t)$ into \eqref{Xeq} gives
\begin{eqnarray}
\left\{
\begin{array}{l}
\displaystyle \bar{A}(t)[\bar{P}(t)Y(t)+p(t)]-F(t)H_3^{-1}(t)F'(t)Y(t)\\
\\
\displaystyle +\bar{B}(t)[I+P(t)\bar{H}(t)]^{-1}[\bar{P}(t)\bar{H}(t)C(t)\bar{P}(t)-\bar{P}(t)\bar{B}'(t)]Y(t)\\
\\
\displaystyle +\bar{B}(t)[I+P(t)\bar{H}(t)]^{-1}[\bar{P}(t)\bar{H}(t)C(t)p(t)+q(t)]\\
\\
\displaystyle=\mathcal{B}_1(t)p(t)+\mathcal{B}_2(t)q(t)\\
\\
\displaystyle +\dot{\bar{P}}(t)Y(t)+\bar{P}(t)[-\bar{A}'(t)Y(t)-C'(t)\bar{H}(t)C(t)(\bar{P}(t)Y(t)+p(t))]\\
\\
\displaystyle +\bar{P}(t)C'(t)\bar{H}(t)[I+P(t)\bar{H}(t)]^{-1}[\bar{P}(t)\bar{H}(t)C(t)\bar{P}(t)-\bar{P}(t)\bar{B}'(t)]Y(t)\\
\\
\displaystyle +\bar{P}(t)C'(t)\bar{H}(t)[I+P(t)\bar{H}(t)]^{-1}[\bar{P}(t)\bar{H}(t)C(t)p(t)+q(t)],
\end{array}\right.\label{XXeq}
\end{eqnarray}
which holds due to our assumptions on $\bar{P}(t)$ and $p(t)$. Substituting \eqref{spt} and \eqref{Z} into the equation for $Y(t)$ in \eqref{hamiltonian2},
shows that it is a linear stochastic differential equation with a unique solution for any $K\in\mathds{R}^n$. This proves the existence of a unique solution triple $(\bar{X}(\cdot),Y(\cdot),\bar{Z}(\cdot))$ of \eqref{hamiltonian2} for any $\xi$. In order to ensure that $\bar{X}(0)=x_0$ for any $x_0\in\mathds{R}^n$, it is necessary and sufficient to have $x_0=\bar{P}(0)K+p(0)$. This equation has a unique solution for any $\xi$ and $x_0$ if and only if $\bar{P}(0)$ is invertible, which we proved is equivalent with the exact controllability of \eqref{main2}.\qed
\section{Minimum-energy control}\label{minsec} In this section we give the solution to the minimum-energy control problem (\ref{ME}). Let us first prove two useful lemmas.
\begin{lemma} Let $(v_1(\cdot),z_1(\cdot))\in\mathcal{A}_{\xi}$ and $(v_2(\cdot),z_2(\cdot))\in\mathcal{A}_{\xi}$ be any two pairs of admissible controls. Then
\begin{eqnarray}
\mathds{E}\int_0^T\Phi(t)F(t)[v_1(t)-v_2(t)]dt=0.\label{lemma1}
\end{eqnarray}
\end{lemma}
{\bf Proof.} By It\^{o}'s product rule, for any admissible pair $(v(\cdot),z(\cdot))\in\mathcal{A}_{\xi}$, we obtain
\begin{eqnarray}
d[\Phi(t)x(t)]&=&[d\Phi(t)]x(t)+\Phi(t)dx(t)-\Phi(t)G(t)[C(t)x(t)+z(t)]dt\nonumber\\
\nonumber\\
&=&-\Phi(t)[A(t)-G(t)C(t)]x(t)dt-\Phi(t)G(t)x(t)dW(t)\nonumber\\
\nonumber\\
&+&\Phi(t)[A(t)x(t)+F(t)v(t)+G(t)z(t)]dt\nonumber\\
\nonumber\\
&+&\Phi(t)[C(t)x(t)+z(t)]dW(t)-\Phi(t)G(t)[C(t)x(t)+z(t)]dt\nonumber\\
\nonumber\\
&=&\!\Phi(t)F(t)v(t)dt\!+\!\Phi(t)\{z(t)\!+\![C(t)\!-\!G(t)]x(t)\}dW(t).\label{Phi}
\end{eqnarray}
Denoting by $x^{(1)}(\cdot)$ and $x^{(2)}(\cdot)$ the solutions to (\ref{main2}) corresponding to $(v_1(\cdot),z_1(\cdot))$ and $(v_2(\cdot),z_2(\cdot))$, respectively, we obtain
\begin{eqnarray}
\Phi(T)\xi-x_0\!\!&=&\!\!\!\int_0^T\!\!\Phi(t)F(t)v_1(t)dt\!+\!\!\int_0^T\!\!\Phi(t)\{z_1(t)\!+\![C(t)\!-\!G(t)]x^{(1)}(t)\}dW(t),\nonumber\\
\nonumber\\
\Phi(T)\xi-x_0\!\!&=&\!\!\!\int_0^T\!\!\Phi(t)F(t)v_2(t)dt\!+\!\!\int_0^T\!\!\Phi(t)\{z_2(t)\!+\![C(t)\!-\!G(t)]x^{(2)}(t)\}dW(t).\nonumber
\end{eqnarray}
The difference of the above two equations is
\begin{eqnarray}
0&=&\int_0^T\Phi(t)F(t)[v_1(t)-v_2(t)]dt\nonumber\\
&+&\int_0^T\Phi(t)\{z_1(t)-z_2(t)+[C(t)-G(t)][x^{(1)}(t)-x^{(2)}(t)]\}dW(t).\nonumber
\end{eqnarray}
The conclusion follows by taking the expectation of both sides.\qed\\

Consider an $\mathds{R}^n$-valued stochastic process $\Gamma(\cdot)$ defined as the solution to the stochastic differential equation
\begin{eqnarray}
\left\{
\begin{array}{l}
\displaystyle d\Gamma(t)=\Gamma_1(t)dt+\Gamma_2(t)dW(t)\\
\\
\displaystyle \Gamma(0)=0,
\end{array}\right.\label{Gamma}
\end{eqnarray}
where $\Gamma_2(\cdot)$ is any process in $L^2_{\mathcal{F}}(0,T;\mathds{R}^n)$, and
\begin{eqnarray}
\Gamma_1(t)\equiv-[\Phi^{-1}(t)]'[C(t)-G(t)]'\Phi'(t)\Gamma_2(t).\nonumber
\end{eqnarray}
\begin{lemma} Let $(v_1(\cdot),z_1(\cdot))\in\mathcal{A}_{\xi}$ and $(v_2(\cdot),z_2(\cdot))\in\mathcal{A}_{\xi}$ be any two pairs of admissible controls. Then
\begin{eqnarray}
\mathds{E}\int_0^T\!\!\Gamma_2'(t)\Phi(t)[z_2(t)-z_1(t)]dt\!=\!-\mathds{E}\int_0^T\!\!\Gamma'(t)\Phi(t)F(t)[v_2(t)-v_1(t)]dt.\label{lemma2}
\end{eqnarray}
\end{lemma}
{\bf Proof.} For any pair of admissible controls $(v(\cdot),z(\cdot))\in\mathcal{A}_{\xi}$,  by It\^{o}'s product rule and (\ref{Phi}), we obtain
\begin{eqnarray}
d[\Gamma'(t)\Phi(t)x(t)]&=&[d\Gamma'(t)]\Phi(t)x(t)+\Gamma'(t)d[\Phi(t)x(t)]+\Gamma_2'(t)\Phi(t)\{z(t)+[C(t)-G(t)]x(t)\}dt\nonumber\\
\nonumber\\
&=&\Gamma_1'(t)\Phi(t)x(t)dt+\Gamma'_2(t)\Phi(t)x(t)dW(t)\nonumber\\
\nonumber\\
&+&\Gamma'(t)\Phi(t)F(t)v(t)dt+\Gamma'(t)\Phi(t)\{z(t)+[C(t)-G(t)]x(t)\}dW(t)\nonumber\\
\nonumber\\
&+&\Gamma_2'(t)\Phi(t)\{z(t)+[C(t)-G(t)]x(t)\}dt\nonumber\\
\nonumber\\
&=&[\Gamma'(t)\Phi(t)F(t)v(t)+\Gamma_2'(t)\Phi(t)z(t)]dt\nonumber\\
\nonumber\\
&+&\Gamma_2'(t)\Phi(t)x(t)dW(t)+\Gamma'(t)\Phi(t)\{z(t)+[C(t)-G(t)]x(t)\}dW(t).\nonumber
\end{eqnarray}
The rest of the proof proceeds as in the previous lemma.\qed
\begin{theorem} (Minimum-energy control) There exists a unique solution to the problem (\ref{ME}) given by
\begin{eqnarray}
v^*(t)&=&H_3^{-1}(t)[F'(t)Y(t)-H_2'Z(t)],\label{optimal v}\\
\nonumber\\
z^*(t)&=&Z(t).\label{optimal z}
\end{eqnarray}
\end{theorem}
{\bf Proof.} We first show that $(v^*(\cdot),z^*(\cdot))\in\mathcal{A}_\xi$. By choosing  the process $\Gamma_2(\cdot)$ in (\ref{Gamma}) as
\begin{eqnarray}
\Gamma_2(t)=[\Phi'(t)]^{-1}[H_1(t)z^*(t)+H_2(t)v^*(t)],\nonumber
\end{eqnarray}
the process $\Gamma_1(\cdot)$ and the equation for $\Gamma(\cdot)$ become
\begin{eqnarray}
\Gamma_1(t)=-[\Phi'(t)]^{-1}[C(t)-G(t)]'[H_1(t)z^*(t)+H_2(t)v^*(t)],\nonumber
\end{eqnarray}
\begin{eqnarray}
\left\{
\begin{array}{l}
\displaystyle d\Gamma(t)=-[\Phi'(t)]^{-1}[C(t)-G(t)]'[H_1(t)z^*(t)+H_2(t)v^*(t)]dt\\
\\
\displaystyle\quad\quad\quad\quad +[\Phi'(t)]^{-1}[H_1(t)z^*(t)+H_2(t)v^*(t)] dW(t)\\
\\
\displaystyle \Gamma(0)=0,
\end{array}\right.\label{Gamma2}
\end{eqnarray}
 Let $\bar{Y}(\cdot)\equiv \Phi'(\cdot)[\Gamma(\cdot)+K]$. We show that $\bar{Y}(\cdot)$ is in fact the process $Y(\cdot)$ of the previous section. Note that $\bar{Y}(0)=\Phi'(0)[\Gamma(0)+K]=K=Y(0)$. The differential of $\bar{Y}(\cdot)$ is:
\begin{eqnarray}
d\bar{Y}(t)&=&[d\Phi'(t)][\Gamma(t)+K]+\Phi'(t)d\Gamma(t)-G'[H_1z^*(t)+H_2(t)v^*(t)]dt\nonumber\\
\nonumber\\
&=&-[A(t)-G(t)C(t)]'\Phi'(t)[\Gamma(t)+K]dt-G'(t)\Phi'(t)[\Gamma(t)+K]dW(t)\nonumber\\
\nonumber\\
&-&[C'(t)-G'(t)][H_1(t)z^*(t)+H_2(t)v^*(t)]dt+[H_1(t)z^*(t)+H_2(t)v^*(t)]dW(t)\nonumber\\
\nonumber\\
&-&G'(t)[H_1(t)z^*(t)+H_2(t)v^*(t)]dt.\nonumber
\end{eqnarray}
After substituting the expressions for $v^*(\cdot)$ and $z^*(\cdot)$, this equation becomes
\begin{eqnarray}
d\bar{Y}(t)&=&\{-[A'(t)-C'(t)G'(t)]\bar{Y}(t)-C'(t)H_2(t)H_3^{-1}(t)F'(t)Y(t)\}dt\nonumber\\
\nonumber\\
&-&C'(t)[H_1(t)-H_2(t)H_3^{-1}H_2'(t)]Z(t)dt\nonumber\\
\nonumber\\
&+&\{-G'(t)\bar{Y}(t)+H_2(t)H_3^{-1}(t)F'(t)Y(t)+[H_1(t)-H_2(t)H_3^{-1}(t)H_2'(t)]Z(t)\}dW(t),\nonumber
\end{eqnarray}
which is satisfied by the process $Y(\cdot)$. Therefore, $\bar{Y}(t)=Y(t)$. Substituting (\ref{optimal v}) and (\ref{optimal z}) in (\ref{main2}) makes it clear that the equations for $x(t)$ and $X(t)$ are the same, and in particular $x(T)=\xi$ a.s.. Hence, $(v^*(\cdot),z^*(\cdot))\in\mathcal{A}_\xi$.\\

We now focus in proving that $(v^*(\cdot),z^*(\cdot))$ are the unique optimal controls. By using equation (\ref{lemma2}), for any control pair $(v(\cdot),z(\cdot))\in\mathcal{A}_{\xi}$, we obtain:
\begin{eqnarray}
&&\mathds{E}\int_0^T\{[H_1(t)z^*(t)+H_2(t)v^*(t)]'[z(t)-z^*(t)]+[H_2'(t)z^*(t)+H_3(t)v^*(t)]'[v(t)-v^*(t)]\}dt\nonumber\\
\nonumber\\
&&=\mathds{E}\int_0^T\{-\Gamma'(t)\Phi(t)F(t)[v(t)-v^*(t)]+[H_2'(t)z^*(t)+H_3(t)v^*(t)]'[v(t)-v^*(t)]\}dt\nonumber\\
\nonumber\\
&&=\mathds{E}\int_0^T[H_2'(t)z^*(t)+H_3(t)v^*(t)-F'(t)\Phi'(t)\Gamma(t)]'[v(t)-v^*(t)]dt\nonumber\\
\nonumber\\
&&=\mathds{E}\int_0^TK'\Phi(t)F(t)[v(t)-v^*(t)]dt=0,\label{mid}
\end{eqnarray}
where the last equality is due to (\ref{lemma1}). For any control pair $(v(\cdot),z(\cdot))\in\mathcal{A}_{\xi}$ we have
\begin{eqnarray}
J(v(\cdot),z(\cdot))&=&\mathds{E}\int_0^T[z'(t)H_1(t)z(t)+2v'(t)H_2'(t)z(t)+v'(t)H_3(t)v(t)]dt\nonumber\\
\nonumber\\
&=&\mathds{E}\int_0^T[z(t)-z^*(t)+z^*(t)]'H_1(t)[z(t)-z^*(t)+z^*(t)]dt\nonumber\\
\nonumber\\
&+&\mathds{E}\int_0^T[v(t)-v^*(t)+v^*(t)]'H_2'(t)[z(t)-z^*(t)+z^*(t)]dt\nonumber\\
\nonumber\\
&+&\mathds{E}\int_0^T[v(t)-v^*(t)+v^*(t)]'H_3(t)[v(t)-v^*(t)+v^*(t)]dt\nonumber
\end{eqnarray}

\begin{eqnarray}
&=&\mathds{E}\int_0^T\{[z(t)-z^*(t)]'H_1(t)[z(t)-z^*(t)]+2[z(t)-z^*(t)]'H_1(t)z^*(t)+(z^*(t))'H_1(t)z^*(t)\}dt\nonumber\\
\nonumber\\
&+&2\mathds{E}\int_0^T\{[v(t)-v^*(t)]'H_2'(t)[z(t)-z^*(t)]+[v(t)-v^*(t)]'H_2'(t)z^*(t)\}dt\nonumber\\
\nonumber\\
&+&2\mathds{E}\int_0^T\{(v^*(t))'H_2'(t)[z(t)-z^*(t)]+(v^*(t))'H_2'(t)z^*(t)\}dt\nonumber\\
\nonumber\\
&+&\mathds{E}\int_0^T\{[v(t)-v^*(t)]'H_3(t)[v(t)-v^*(t)]+2(v^*(t))'H_3(t)[v(t)-v^*(t)]+(v^*(t))'H_3(t)v^*(t)\}dt\nonumber\\
\nonumber\\
&=&J(v^*(\cdot),z^*(\cdot))+\mathds{E}\int_0^T\left[
\begin{array}{l}
z(t)-z^*(t)\\
v(t)-v^*(t)
\end{array}\right]'\left[
\begin{array}{ll}
H_1(t)&H_2(t)\\
H_2'(t)&H_3(t)
\end{array}\right]\left[
\begin{array}{l}
z(t)-z^*(t)\\
v(t)-v^*(t)
\end{array}\right]dt\nonumber\\
\nonumber\\
&+&2\mathds{E}\int_0^T\{[z(t)-z^*(t)]'[H_1(t)z^*(t)+H_2(t)v^*(t)]+[v(t)-v^*(t)]'[H_2'(t)z^*(t)+H_3(t)v^*(t)]\}dt,\nonumber
\end{eqnarray}
Due to (\ref{mid}), for any $(v(\cdot),z(\cdot))\in\mathcal{A}_{\xi}$ we have
\begin{eqnarray}
J(v(\cdot),z(\cdot))&=&J(v^*(\cdot),z^*(\cdot))+\mathds{E}\int_0^T\left[
\begin{array}{l}
z(t)-z^*(t)\\
v(t)-v^*(t)
\end{array}\right]'\left[
\begin{array}{ll}
H_1(t)&H_2(t)\\
H_2'(t)&H_3(t)
\end{array}\right]\left[
\begin{array}{l}
z(t)-z^*(t)\\
v(t)-v^*(t)
\end{array}\right]dt\nonumber\\
&\geq&J(v^*(\cdot),z^*(\cdot)),\nonumber
\end{eqnarray}
with equality if and only if $v(t)=v^*(t)$, $a.e.$ $t\in[0,T]$ $a.s.$, and $z(t)=z^*(t)$, $a.e.$ $t\in[0,T]$ $a.s.$.\qed
\section{Stochastic LQ regulator with a fixed final state}
Ever since its introduction by Kalman~\cite{kalman1}, the LQ regulator has been studied extensively in both deterministic~\cite{anderson},~\cite{athans}, and stochastic~\cite{won1},~\cite{won2},~\cite{yong}, settings. The version of this regulator with a {\it fixed final state}~\cite{lewis} is the minimum-energy control problem with the cost functional that has a penalty on the state as well as on the control. As an extension of our results on minimum-energy control, we give the solution to the following {\it stochastic} LQ regulator problem with a fixed final state.\\

{\bf LQ regulator with a fixed final state.} {\it Let $Q(\cdot)\in L^\infty(0,T;\mathds{R}^{n\times n})$ be a given symmetric matrix such that $Q(t)\geq0$, $a.e.$ $t\in[0,T]$. For any given $x_0\in\mathds{R}^n$ and $\xi\in L^2(\Omega,\mathcal{F}_T,P;\mathds{R}^n)$ find the control process $u(\cdot)\in\mathcal{U}_{\xi}$ that minimizes the cost functional}
\begin{eqnarray}
\widehat{J}(u(\cdot))=\mathds{E}\int_0^T[x'(t)Q(t)x(t)+u'(t)R(t)u(t)]dt.\label{main costlq}
\end{eqnarray}

We solve this problem by transforming it into an equivalent minimum-energy control problem of the previous sections. Consider the Riccati differential  equation
\begin{eqnarray}
\left\{
\begin{array}{l}
\displaystyle \dot{{P}}(t)+{P}(t)A(t)+A'(t){P}(t)+C'(t)P(t)C(t)+Q(t)\\
\\
\displaystyle -[{P}(t)B(t)+C'(t)P(t)D(t)][D'(t)P(t)D(t)+R(t)]^{-1}[B'(t)P(t)+D'(t)P(t)C(t)]=0,\\
\\
\displaystyle {P}(T)=0,\\
\\
\displaystyle D'(t)P(t)D(t)+R(t)>0,
\end{array}\right.\nonumber
\end{eqnarray}
which has a unique solution (see, e.g.,~\cite{won1},~\cite{won2},~\cite{yong},~\cite{XYZ}). We introduce the following matrices for notational convenience:
\begin{eqnarray}
&&\widehat{A}(t)\equiv A(t)-B(t)[D'(t)P(t)D(t)+R(t)]^{-1}[B'(t)P(t)+D'(t)P(t)C(t)],\nonumber\\
\nonumber\\
&&\widehat{B}(t)\equiv B(t),\nonumber\\
\nonumber\\
&&\widehat{C}(t)\equiv C(t)-D(t)[D'(t)P(t)D(t)+R(t)]^{-1}[B'(t)P(t)+D'(t)P(t)C(t)],\nonumber\\
\nonumber\\
&&\widehat{D}(t)\equiv D(t),\nonumber\\
\nonumber\\
&&\widehat{R}(t)\equiv D'(t)P(t)D(t)+R(t).\nonumber
\end{eqnarray}

By introducing a new control $\widehat{u}(t)\equiv u(t)+[D'(t)P(t)D(t)+R(t)]^{-1}[B'(t)P(t)+D'(t)P(t)C(t)]x(t)$, we can rewrite the state equation (\ref{main system}) as
\begin{eqnarray}
\left\{
\begin{array}{l}
\displaystyle dx(t)=[\widehat{A}(t)x(t)+\widehat{B}(t)\widehat{u}(t)]dt+[\widehat{C}(t)x(t)+\widehat{D}(t)\widehat{u}(t)]dW(t)\\
\\
\displaystyle x(0)=x_0\in\mathds{R}^n,\quad\mbox{is given}.
\end{array}\right.\label{mainlq}
\end{eqnarray}
Using the completion of squares method of stochastic LQ control~\cite{yong}, we can rewrite the cost functional (\ref{main costlq}) as
\begin{eqnarray}
\widehat{J}(\widehat{u}(\cdot))=x_0'P(0)x_0+\mathds{E}\int_0^T\widehat{u}'(t)\widehat{R}(t)\widehat{u}(t)dt.\label{trans}
\end{eqnarray}
Apart from the obvious change of notation, the problem of minimizing (\ref{trans}) subject to (\ref{mainlq}) and $\widehat{u}(\cdot)\in\mathcal{A}_\xi$, is the stochastic minimum-energy control problem. Hence, its solution can be obtained by applying the results of the previous sections.
\section{Conclusions}
A stochastic version of the classical minimum-energy control problem is formulated. The system is driven by a Brownian motion and all coefficients can be time-varying. By assuming the exact controllability of the system, complete solution is given. The minimum-energy control problem is crucial is solving several optimal control problems involving terminal state constraints and quadratic criteria. One such problem is the stochastic LQ regulator with a fixed final state, and it has been solved as an application of our results on minimum-energy control. We expect that the method proposed in this paper will prove useful in tackling more general minimum-energy control problems.

\end{document}